\documentclass[11pt]{article}
\usepackage[T1]{fontenc}
\usepackage[utf8]{inputenc}
\usepackage{lmodern,microtype}
\usepackage[a4paper,margin=30mm]{geometry}
\usepackage{amsmath,amssymb,amsthm,mathtools,mathrsfs,array,booktabs,needspace}
\usepackage{aliascnt}
\usepackage{xcolor,hyperref}
\usepackage[nameinlink,capitalise,noabbrev]{cleveref}
\usepackage{fancyhdr}
\hypersetup{colorlinks=true,linkcolor=blue!45!black,citecolor=green!35!black,
 urlcolor=blue!55!black,pdftitle={The Exact Nef-Curve Cone Formula for dlt Pairs},
 pdfsubject={Exact rational generators for nef-curve cones of dlt pairs}}
\newtheorem{theorem}{Theorem}[section]
\newaliascnt{question}{theorem}\newtheorem{question}[question]{Question}\aliascntresetthe{question}
\crefname{question}{Question}{Questions}\Crefname{question}{Question}{Questions}
\newaliascnt{proposition}{theorem}\newtheorem{proposition}[proposition]{Proposition}\aliascntresetthe{proposition}
\newaliascnt{lemma}{theorem}\newtheorem{lemma}[lemma]{Lemma}\aliascntresetthe{lemma}
\newaliascnt{corollary}{theorem}\aliascntresetthe{corollary}
\theoremstyle{remark}
\newaliascnt{remark}{theorem}\aliascntresetthe{remark}
\newcommand{\C}{\mathbb C}\newcommand{\Q}{\mathbb Q}\newcommand{\R}{\mathbb R}
\newcommand{\PP}{\mathbb P}\newcommand{\cO}{\mathcal O}
\newcommand{\NE}{\overline{\operatorname{NE}}_1}\newcommand{\NM}{\overline{\operatorname{NM}}_1}\newcommand{\Mov}{\operatorname{Mov}}

\newcommand{\reg}{\mathrm{reg}}

\title{\bfseries The Exact Nef-Curve Cone Formula for dlt Pairs}
\author{Jingcao Wu\thanks{School of Mathematics, Shanghai University of Finance and Economics, Shanghai 200433, People's Republic of China. E-mail: \texttt{wujincao@shufe.edu.cn}.}}
\date{}
\begin{document}\maketitle
\begin{abstract}
We prove an exact rational-curve formula for the nef-curve cone of every complex projective dlt pair with effective real boundary, answering the related question proposed by Cascini. Every $(K_X+\Delta)$-negative extremal ray is generated by a movable rational curve whose normalization is free in the smooth locus. We reduce the klt case to known results and treat the remaining limiting rays with a new argument based on an orthogonal movable divisor and induction on dimension.
\end{abstract}
\medskip\noindent\textbf{2020 Mathematics Subject Classification.}
Primary 14E30; Secondary 14J17, 14C20, 14J45.
\par\noindent\textbf{Keywords.}
Nef-curve cone; movable rational curves; dlt pair; free rational curves.

\section{Introduction}\label{sec:intro}

The cone theorem for nef curves describes the outcome of the minimal model program in terms of numerical classes of movable curves. For a normal projective variety $X$, denote by $\NE(X)$ the closed cone of effective curves in $N_1(X)$ and by $\NM(X)$ the closed cone of nef curves, dual to the pseudoeffective cone of Cartier divisors. For a numerical divisor $L$, write $\NE(X)_{L\ge0}$ for the intersection with the indicated closed half-space. The results of Batyrev, Araujo and Lehmann \cite{Araujo09,Batyrev92,Lehmann12} motivate the following question, which is the starting point of this paper.

\begin{question}\label{que:cascini}
Let $(X,\Delta)$ be a projective divisorially log terminal (dlt) pair over $\C$. Do there exist at most countably many movable rational curves $C_i$ with $(K_X+\Delta)\cdot C_i<0$ such that
\begin{equation}\label{eq:cascini}
\NE(X)_{K_X+\Delta\ge0}+\NM(X)=\NE(X)_{K_X+\Delta\ge0}+\sum_i\R_{\ge0}[C_i]?
\end{equation}
Here the sum on the right denotes finite nonnegative linear combinations; no closure is taken.
\end{question}

Cascini \cite[Theorem~2.23]{Cascini13} collects a dlt cone theorem from \cite{Araujo09,Batyrev92,Lehmann12} with countably many movable-curve generators. The discussion following that theorem identifies two further assertions in \eqref{eq:cascini}---that the generators can be chosen rational and that the equality holds without taking a closure---as conjectural; at the time, they were known in
dimensions at most three \cite[Theorem~2.23 and the discussion following
it]{Cascini13}. Lehmann's dlt cone theorem also uses movable curves without asserting that they are rational \cite[Theorem~1.3]{Lehmann12}. For klt pairs, Das proves the exact equality without closure, with movable-curve generators \cite[Theorem~1.9]{Das20}.

The remaining arbitrary-dimensional case is the representation of limiting negative rays for dlt pairs: klt perturbations yield only a closed sum in the limit and do not by themselves produce rational curves on those rays; see \cref{subsec:known-consequences}. Restricting the generators to negative degree causes no loss, since every movable curve of nonnegative $(K_X+\Delta)$-degree already belongs to the first term on the right of
\eqref{eq:cascini}. We impose no restriction on dimension or Picard number.

\begin{theorem}[Exact rational generators for dlt pairs]\label{thm:exact-dlt}
Let $(X,\Delta)$ be a complex projective dlt pair with effective real boundary such that $K_X+\Delta$ is $\mathbb R$-Cartier, and put $L=K_X+\Delta$. There are at most countably many movable rational curves $C_i$ with $L\cdot C_i<0$ such that
\[
\NE(X)_{L\ge0}+\NM(X)=\NE(X)_{L\ge0}+\sum_i\mathbb R_{\ge0}[C_i],
\]
where the sum consists of finite nonnegative linear combinations. Moreover, every $L$-negative extremal ray of $\NE(X)_{L\ge0}+\NM(X)$ is generated by a movable rational curve whose normalization is free in $X_{\rm sm}$.
\end{theorem}

The theorem answers \cref{que:cascini} affirmatively in every dimension, with equality as an algebraic sum. The known klt cone formula \cite[Theorem~1.9]{Das20} and free rational curves on klt Fano fibres \cite[Theorem~1.1]{Xu26} handle the klt rays. The new step is to represent limiting negative rays that are not preserved by klt perturbation. We construct a movable divisor orthogonal to each such ray, then use the pseudo-effective reduction map, induction on dimension, and dlt descent to obtain its rational representative.

The proof has two stages. \cref{sec:cone} records the klt MMP input and the dlt perturbation formula. \cref{sec:exact} treats the limiting rays, first on a log-smooth model and then by descent to an arbitrary dlt pair. Standard deformation results for free maps provide the freeness and codimension-two avoidance needed in both stages.

\section{Mori fibre spaces and the exact nef-curve cone}\label{sec:cone}

We establish the klt input needed for the treatment of limiting rays in \cref{sec:exact}.

\subsection{Rational representatives of MMP rays}\label{subsec:mmp-rays}

For a projective normal variety $X$, write $N^1(X)$ and $N_1(X)$ for its real numerical divisor and curve spaces. We use the usual rational structures given by Cartier divisors and numerical curve classes. For a numerical divisor class $L$, put
\[
B_L(X)=\NE(X)_{L\ge0},\qquad\mathcal C_L(X)=B_L(X)+\NM(X).
\]
The cone $\NM(X)$ is dual to the pseudoeffective cone of Cartier divisors and is contained in $\NE(X)$. A movable curve is an irreducible curve in an algebraic family whose evaluation map dominates $X$; its class lies in $\NM(X)$. Here, rationality refers to the normalization of the image curve, not merely to its numerical class. A nonconstant map from $\PP^1$ pushes forward to a positive integral multiple of its reduced rational image, and the two classes span the same ray.

Fix an ample Cartier divisor $H$ on $X$. The slice
\[
 S_H=\{\alpha\in\NE(X):H\cdot\alpha=1\}
\]
is compact by Kleiman's criterion. Both summands of $\mathcal C_L(X)$ are closed cones in $\NE(X)$; their sum is closed, since a sequence of bounded sums has bounded summands after applying $H$. We normalize every nonzero ray in this cone by its unique point in $S_H$. In particular, countably many rational rays suffice whenever a cone is proved to be generated by a collection of rational rays: the numerical lattice has only countably many such rays.

\paragraph{Free maps and their deformations.}
All varieties and stacks are over $\C$. A nonconstant map from the ordinary source $\PP^1$ to a smooth variety or Deligne--Mumford stack is \emph{free} if $f^*T_X$ is globally generated, equivalently
\[
 H^1\bigl(\PP^1,f^*T_X(-1)\bigr)=0.
\]
Freeness does not require an embedding.

For a dominating family of ordinary rational maps to a smooth Deligne--Mumford stack, generic smoothness of evaluation gives sections of the pullback tangent bundle that span its generic fibre. A general map is therefore free. For a free map, the full Hom space and its one-point evaluation map are smooth, since $H^1(\PP^1,f^*T_X(-p))=0$ for every
point $p$. The incidence codimension argument then shows that general deformations in the full Hom space avoid any fixed closed subset of codimension at least two. If a free map is not contained in a smooth Cartier divisor, general deformations meet that divisor transversely: the finite incidence projection is generically smooth in characteristic zero. The same statements hold on the regular locus of a singular stack. They apply to the full free deformation space, not to an arbitrary restricted subfamily (cf. \cite{AGV08,Kollar96}).

\paragraph{The MMP pullback.}
Let $(X,\Gamma)$ be projective $\Q$-factorial klt, and consider an MMP with scaling for $K_X+\Gamma$ ending in a Mori fibre space
\[
\phi:X\dashrightarrow X'\xrightarrow{p}Z,\qquad \Gamma'=\phi_*\Gamma.
\]
We use the numerical pullback $\phi^*_{\rm num}:N_1(X')\to N_1(X)$ dual to pushforward of divisor classes, as in \cite[Section~4]{Araujo09}. We only apply it to a curve on a general fibre of $p$. All nonzero effective curves in fibres have the same relative numerical direction because $\rho(X'/Z)=1$. The ray of their numerical pullbacks is well-defined.

We now show that this ray has a movable rational representative. The effective divisor $\Gamma'$ has nonnegative degree on a movable curve in a general fibre. Since $\rho(X'/Z)=1$ and $-(K_{X'}+\Gamma')$ is $p$-ample, $-K_{X'}$ is $p$-ample as well. By generic restriction and adjunction over the smooth locus of the base, the general positive-dimensional fibre $F$ is a normal klt Fano variety. Since $X'$ is $\Q$-factorial, $K_{X'}$ is also $\Q$-Cartier. Xu's theorem \cite[Theorem~1.1]{Xu26} gives a nonconstant free map $f:\PP^1\to F_{\reg}$. After restricting to the smooth locus of the base and the flat locus of $p$, a point regular on its fibre is a smooth point of the morphism. Thus $p$ is smooth along $f$, and the tangent sequence along $f$ is
\[
0\longrightarrow f^*T_{F_{\reg}}\longrightarrow f^*T_{X'_{\reg}}\longrightarrow \cO_{\PP^1}^{\oplus\dim Z}\longrightarrow0.
\]
The first term is globally generated, so the extension splits on $\PP^1$ and $f$ is free in $X'_{\reg}$. In the full free component, the curve moves dominantly and remains vertical: its degree against the pullback of an ample divisor on $Z$ is zero. The inverse of an MMP birational contraction is an isomorphism to an open subset of $X$ away from a closed subset of codimension at least two in $X'$. A general free deformation avoids that subset by the deformation facts above. The curve therefore transfers to a free rational map in $X_{\reg}$. For every Cartier divisor $D$ on $X$, the transferred curve has the same degree against $D$ as the curve in $X'$ has against $\phi_*D$. Thus its class is the numerical pullback, not merely a limit of classes in that direction. In particular,
\begin{equation}\label{eq:transfer-degree}
-(K_X+\Gamma)\cdot C=-(K_{X'}+\Gamma')\cdot C'>0
\end{equation}
for the corresponding map cycles $C$ and $C'$; the same equality holds for their reduced images. This realizes the rational-curve input identified in \cite[Remark~3.4]{Lehmann12}.

When $X$ is not $\Q$-factorial but $(X,\Gamma)$ is klt, take a projective small $\Q$-factorialization $q:\widehat X\to X$ by \cite{BCHM10}. For the strict transform $\widehat\Gamma$ one has $K_{\widehat X}+\widehat\Gamma=q^*(K_X+\Gamma)$, and the pair is klt. We use the pushforward identities
\[
q_*\NE(\widehat X)=\NE(X),\qquad q_*\NM(\widehat X)=\NM(X)
\]
from \cite[Lemma~2.1]{Lehmann12}. The analogous identity for the nonnegative effective cone holds for a pulled-back $L$, since the sign is unchanged on every lift of a numerical class. As $q$ is small, free rational representatives on $\widehat X$ can be deformed away from its exceptional locus and transferred to $X_{\reg}$, with their pushforward classes unchanged. Thus the $\Q$-factorial arguments below descend to the
original cone.

\subsection{Klt inputs and the dlt perturbation gap}
\label{subsec:known-consequences}

For a fixed projective $\Q$-factorial klt pair $(X,\Gamma)$ put $L=K_X+\Gamma$. Denote by $\mathscr R(L)$ the rays obtained by numerically pulling back fibre rays of $(K_X+\Gamma)$-MMPs that end in Mori fibre spaces. Araujo's cone theorem \cite[Theorem~1.1]{Araujo09} gives
\begin{equation}\label{eq:araujo-closed}
 \mathcal C_L(X)=B_L(X)+\overline{\sum_{R\in\mathscr R(L)} R}.
\end{equation}
If $L$ is pseudoeffective, no negative rays are needed. Otherwise, every ray in $\mathscr R(L)$ is $L$-negative and has a free rational representative by the Mori fibre argument in \cref{subsec:mmp-rays}.

We need to identify the negative extremal directions in \eqref{eq:araujo-closed} with actual MMP rays. When $(X,\Gamma)$ is $\epsilon$-lc, Das \cite[Proposition~7.1]{Das20} gives a constant $G=G(\dim X,\epsilon)$ such that each $R\in\mathscr R(L)$ has a movable integral curve representative $C$ with $-L\cdot C\le G$. If $u_R$ is the point of $R$ normalized by $H\cdot u_R=1$ and $-L\cdot u_R\ge\eta>0$, then
\begin{equation}\label{eq:das-bound}
 H\cdot C\le G/\eta.
\end{equation}
There are only finitely many integral numerical curve classes of bounded $H$-degree. Indeed, a rational Cartier basis $D_i$ of $N^1(X)$ and integers $m_i$ with $m_iH\pm D_i$ ample bound the integer pairings $D_i\cdot C$. Thus the normalized MMP rays are locally finite away from the nonnegative face. On the compact $H$-slice, Araujo's formula places
each negative extreme point in the closure of the normalized MMP-ray set; local finiteness then makes it an actual MMP ray. Nonnegative extreme points lie in $B_L(X)$. Minkowski's and Carath\'eodory's theorems give the exact finite-sum formula using $B_L(X)$ and the actual MMP rays. The bound $G$ also supplies the uniform curve-degree estimate
used in \cref{prop:lift}.

For dlt pairs, Lehmann's klt perturbation \cite[Lemma~4.6]{Lehmann12} gives the known closure-level rational generator statement described in the introduction. The remaining issue is to represent the limiting negative rays without taking a closure.

\section{Proof of the exact dlt cone formula}\label{sec:exact}
It suffices to represent each negative extremal ray by a movable rational curve. Indeed, the slice $\mathcal C_L(X)\cap S_H$ is compact; its nonnegative extreme rays lie in $B_L(X)$, while each negative extreme ray lies in $\NM(X)$. The compact-slice argument in
\cref{subsec:mmp-rays} then gives the finite-sum formula and countability asserted in \cref{thm:exact-dlt}.

We first prove the ray statement for log-smooth pairs, distinguishing rays that persist under small perturbations from limiting rays. The exact klt formula and rational representatives of MMP rays were established in \cref{subsec:mmp-rays,subsec:known-consequences}.

Fix a negative extremal ray $R=\mathbb R_{\ge0}\alpha$ and normalize $H\cdot\alpha=1$. Put
\[
T_\alpha=\{v\in\NE(X): H\cdot v=1,\ \alpha-tv\in\NM(X)\text{ for some }t>0\},\qquad \eta=\inf_{v\in T_\alpha}(-L\cdot v).
\]
For every numerical divisor $L'$ with $L'\cdot\alpha<0$, the ray persists as an extremal ray of $B_{L'}(X)+\NM(X)$ exactly when
\begin{equation}\label{eq:test}
R\text{ is extremal in }\mathcal C_{L'}(X)\quad\Longleftrightarrow\quad L'\cdot v<0\quad\text{for every }v\in T_\alpha .
\end{equation}
Indeed, if the inequality fails, $\alpha=(\alpha-tv)+tv$ splits into a nonzero movable class and a nonzero class in $B_{L'}$ outside $R$. Conversely, in any decomposition of $\alpha$ in $B_{L'}+\NM$, a nonzero $B_{L'}$ summand, normalized by its $H$-degree, belongs to $T_\alpha$ and has nonnegative $L'$-degree; if there is no such summand, extremality in $\NM$ applies. In particular $L\cdot v<0$ throughout $T_\alpha$, so $\eta\ge0$. We use the closure of $T_\alpha$ only when taking limits.

We work first with a smooth projective variety and a simple normal crossings boundary. We apply the following dichotomy to each fixed negative extremal ray.

\begin{lemma}\label{lem:div}
If $\eta>0$, the given ray persists for all sufficiently small numerical perturbations of $L$ and has a movable rational representative. If $\eta=0$, there is $v\in\overline{T_\alpha}$ with $L\cdot v=0$ and a nonzero $D\in\overline{\Mov}^{1}(X)$ satisfying $D\cdot\alpha=D\cdot v=0$.
\end{lemma}
\begin{proof}
For $\eta>0$, it suffices in \eqref{eq:test} to choose $L'$ so that
\[
\sup_{z\in\NE(X),\,H\cdot z=1}|(L'-L)\cdot z|<\eta.
\]
Choose such an $L'=K_X+\Gamma$ by approximating the coefficients of the simple normal crossings boundary $\Delta$ by nonnegative rational numbers strictly below one. Then $(X,\Gamma)$ is klt, and $L'$ is not pseudoeffective since $L'\cdot\alpha<0$. The exact klt formula expresses $\alpha$ as a sum of a class in the nonnegative cone and classes on actual Mori fibre rays. Extremality forces one of these rays to be $R$, giving the
required rational representative.

Suppose $\eta=0$. Compactness gives $v\in\overline{T_\alpha}$ with $L\cdot v=0$, and \eqref{eq:test} shows $v\notin T_\alpha$. Every pseudo-effective divisor $B$ with $B\cdot\alpha=0$ satisfies $B\cdot w\le0$ for $w\in T_\alpha$, hence $B\cdot v\le0$. For each $k\ge1$ the class $\alpha-v/k$ lies outside $\Mov_1(X)$. Since this cone is dual to the closure of the cone generated by effective divisors, there is a prime divisor $E_k$ with
\[
     0< E_k\cdot\alpha < E_k\cdot v/k .
\]
Here $E_k\cdot\alpha\ge0$ by movability, and equality would force $E_k\cdot v\le0$ by the preceding observation. No fixed prime can satisfy this strict inequality for arbitrarily large $k$. Pass to distinct primes and normalize their classes in a compact base of the pseudo-effective cone, cut out by a linear functional strictly positive on nonzero pseudo-effective classes. A subsequential limit $D$ is nonzero and satisfies $D\cdot\alpha=0$. The same inequalities give $D\cdot v\ge0$, while $D\cdot\alpha=0$ gives $D\cdot v\le0$.

It remains to show that $D$ is movable. Let $M$ be a rational Cartier divisor whose class lies in $D+\operatorname{Amp}(X)$. Choose two distinct primes $E_i,E_j$ sufficiently far out in the subsequence, and positive rational numbers $a_i,a_j$ such that $M-a_iE_i$ and $M-a_jE_j$ are ample. Such scalars exist by openness, since the normalized classes tend to $D$. For sufficiently divisible $p$, the subsystem $pa_iE_i+|p(M-a_iE_i)|$ has no fixed prime outside $E_i$; the analogous subsystem for $E_j$ has no fixed prime outside $E_j$. Hence the full system $|pM|$ has no fixed prime. Letting the class of $M$ approach $D$ proves that $D$ lies in the closed movable divisor cone.
\end{proof}

When $\eta=0$, an ample perturbation that keeps $R$ negative destroys its extremality. Indeed, for every ample $A$ with $(L+A)\cdot\alpha<0$, the minimum of $A$ on the $H$-unit slice is positive. Some $v'\in T_\alpha$ therefore satisfies $(L+A)\cdot v'>0$, so
\eqref{eq:test} rules out extremality for $L+A$. Certain boundary-specific perturbations can still preserve $R$. For example, if every component of $\Delta$ with irrational coefficient is nef, rounding those coefficients down to rational values gives a rational boundary $\Gamma$ with $\Delta-\Gamma$ nef. Then $B_{K_X+\Gamma}\subset B_L$, and $R$ remains negative and extremal in the smaller cone. The lemma handles the zero-gap case without requiring a boundary perturbation that preserves $R$.

\subsection{Movable divisors on a sequence of birational models}

Fix a smooth projective $X$, an ample Cartier divisor $H$, and a real class $D\in\overline{\Mov}^1(X)$. For a prime divisor $E$ on a smooth birational model, write $\sigma_E$ for the asymptotic order of vanishing along $E$. For a big rational divisor, it is the infimum of the vanishing orders of its complete linear systems, normalized by their multiples. It extends to the pseudoeffective boundary by ample perturbation. We write $P_\sigma(M)=M-N_\sigma(M)$ for the positive part of the divisorial Zariski decomposition \cite[Definitions~2.9--2.13]{LehmannReduction14}.

\begin{lemma}\label{lem:approx}
Given a smooth projective birational model $\pi_0:W_0\to X$, there is a tower $\pi_j:W_j\to X$ of smooth models, with morphisms $h_{jr}:W_j\to W_r$ for $j\ge r$, and semiample rational divisors $P_j$ on $W_j$ for $j\ge1$, with the following properties:
\begin{enumerate}
\item For every fixed $r\ge0$, $(h_{jr})_*P_j\to P_\sigma(\pi_r^*D)$ numerically as $j\to\infty$.
\item If $\alpha\in\Mov_1(X)$, $D\cdot\alpha=0$ and $H\cdot\alpha=1$, then, setting $h_{mm}=\mathrm{id}$,
\[
S_m=\sum_{j=1}^m h_{mj}^*P_j \quad\text{satisfies}\quad 0\le S_m\cdot\pi_m^!\alpha\le 1.
\]
\end{enumerate}
The models may simultaneously resolve a fixed boundary together with the exceptional divisors.
\end{lemma}
\begin{proof}
Put $q_j=2^{-j}$. Choose rational Cartier divisors $B_j$ with
\[
 D+\tfrac12q_jH <_{\rm amp} [B_j] <_{\rm amp} D+q_jH,
\]
where $U<_{\rm amp}V$ means that $V-U$ is ample. The open region is
nonempty, and rational numerical classes are dense, so each $B_j$ can be
chosen big. For sufficiently divisible $k$, the system $|kB_j|$ has no
fixed prime divisor. Indeed, approximate $D$ by a rational nonnegative
sum $V_j$ of divisors whose linear systems have no fixed prime, closely
enough that $B_j-V_j$ is ample. This is possible by the definition of the
closed movable cone and the openness of ampleness. After clearing
denominators, products of sections of the moving summands and the ample
summand give the required subsystem.

Suppose $W_{j-1}$ has been chosen. For its finitely many
$\pi_{j-1}$-exceptional primes $E$, choose a common divisible multiple
$k_j$ such that $|k_jB_j|$ has no fixed prime on $X$ and
\[
0\le k_j^{-1}\operatorname{ord}_E\bigl(|k_j B_j|\bigr)-\sigma_E(\pi_{j-1}^*B_j) < 1/j.
\]
Here the order is the minimum vanishing order of the pulled-back system.
The infimum formula gives a suitable multiple for each $E$; passing to a
common multiple preserves the bounds. Resolve the base ideal and the
boundary on a smooth model $W_j$ dominating $W_{j-1}$. Write
$P_j=\pi_j^*B_j-F_j/k_j$, where $F_j$ is the fixed divisor of the
pulled-back system. Then $k_jP_j$ is basepoint free, and $F_j$ is
exceptional over $X$.

Fix $r$. Every $\pi_r$-exceptional prime persists by strict transform on $W_{j-1}$ for $j>r$. For each such prime $E$, lower semicontinuity gives $\sigma_E(\pi_r^*D)\le\liminf_j\sigma_E(\pi_r^*B_j)$. The reverse upper bound follows from subadditivity since $\pi_r^*([B_j]-D)$ is nef and has zero asymptotic order. The orders for $B_j$ can be computed on either model since their complete systems pull back. It follows that the coefficient of $(h_{jr})_*F_j/k_j$ at $E$ tends to $\sigma_E(\pi_r^*D)$. All other coefficients of this pushforward are zero. Also $N_\sigma(\pi_r^*D)$ is exceptional over $X$, since $D$ is movable in codimension one. This proves the first assertion. Finally $(\pi_j)_*P_j=B_j$ numerically, and
\[
0\le S_m\cdot\pi_m^!\alpha=\sum_{j=1}^m B_j\cdot\alpha\le\sum_{j=1}^m q_j\le1.
\]
\end{proof}

Choose $W_0$ so that the pseudo-effective reduction of $D$ is represented by a morphism $f:W_0\to Z$ with connected fibres and normal projective base. On $W_m$ call a curve horizontal if $f_* (h_{m0})_*[C]\ne0$.
\begin{lemma}\label{lem:vertical}
For $b>0$, let $d_m(b)$ be the infimum of $S_m\cdot C$ over horizontal integral curves in covering families on $W_m$ with
\[
  H\cdot (\pi_m)_* C\le b,\qquad K_{W_m}\cdot C\le0.
\]
With $\inf\varnothing=+\infty$, one has $d_m(b)\to+\infty$ for every fixed $b$.
\end{lemma}
\begin{proof}
Otherwise there is an unbounded sequence of indices $k$ and such curves $C_k$ with $S_k\cdot C_k\le A$ for a fixed $A$. Bounded ample degree gives only finitely many integral numerical curve classes on $X$ in the effective cone. Pass to a subsequence with $(\pi_k)_*[C_k]=c$ constant. Since \(C_k\) is movable and \(K_{W_k/X}\) is effective,
\[
0\le K_{W_k/X}\cdot C_k= K_{W_k}\cdot C_k-K_X\cdot c\le -K_X\cdot c=:b'.
\]
For each fixed \(j\), the classes \(c_{k,j}=(h_{kj})_*[C_k]\) are represented by members of covering families on \(W_j\). A general member of the original family meets the locus on which \(h_{kj}\) is an isomorphism, so its image has multiplicity one. The numerical class stays constant in the family. Each \(c_{k,j}\) pairs nonnegatively with the prime exceptional divisors of \(\pi_j\), all of which occur in \(K_{W_j/X}\) with positive integer coefficients. The identity
\[
 K_{W_k/X}=h_{kj}^*K_{W_j/X}+K_{W_k/W_j}
\]
bounds each of these integral pairings by \(b'\). They and \(c\) determine \(c_{k,j}\), since \(N^1(W_j)_{\R}\) is spanned by \(\pi_j^*N^1(X)_{\R}\) and the exceptional prime divisors. There are therefore finitely many choices on each fixed \(W_j\). A diagonal subsequence makes \(c_{k,j}\) eventually constant for every \(j\), of value \(c_j\). These classes are compatible under pushforward, and \(f_*c_0\ne0\).

The bounded nondecreasing integers \(K_{W_j/X}\cdot c_j\) eventually stabilize, say for \(j\ge r\). Thus \(K_{W_j/W_r}\cdot c_j=0\) for all \(j\ge r\). Each \(c_j\) annihilates every prime exceptional divisor of \(h_{jr}\), by movability and the strict positivity of
its coefficient in \(K_{W_j/W_r}\). Equivalently, \(c_j\) is the numerical pullback of \(c_r\).

For each fixed $J$ the classes $(h_{kj})_*[C_k]$ equal $c_j$ for all $1\le j\le J$ simultaneously once $k$ is sufficiently large in the subsequence. Every summand of $S_k$ is nef, hence
\[
  0\le\sum_{j=1}^J P_j\cdot c_j\le S_k\cdot C_k\le A.
\]
The nonnegative series $\sum_j P_j\cdot c_j$ thus converges, so $P_j\cdot c_j\to0$. For $j\ge r$, numerical pullback gives $P_j\cdot c_j=((h_{jr})_*P_j)\cdot c_r$. \cref{lem:approx} on the fixed model $W_r$ yields
\[
      P_\sigma(\pi_r^*D)\cdot c_r=0.
\]

Choose one of the projected covering families on \(W_r\) realizing \(c_r\), and let its
very general image curve move in the induced covering family on \(X\). Its strict transform has class \(c_r\). The implication \((2)\Rightarrow(1)\) of \cite[Theorem~5.2]{LehmannReduction14} says that the restricted numerical dimension of \(D\) along this very general curve is zero. By the description of the pseudoeffective reduction
map as the generic quotient by precisely these families \cite[Theorem~1.5]{LehmannReduction14}, the strict transforms on \(W_0\) are contracted by \(f\). Their numerical class is \(c_0\), contrary to \(f_*c_0\ne0\).
\end{proof}

\subsection{Lifting the ray to the divisor fibre}

\begin{proposition}\label{prop:lift}
Let $(X,\Delta)$ be a projective log smooth pair with boundary coefficients at most one, put $L=K_X+\Delta$, and let $\alpha$ span an $L$-negative extremal ray of $\mathcal C_L(X)$. Normalize by $H\cdot\alpha=1$ for an ample Cartier divisor $H$. Suppose that $D$ is a nonzero class in the closed movable divisor cone with $D\cdot\alpha=0$. Let $\pi_0:W_0\to X$ be a log resolution carrying a morphism $f:W_0\to Z$ representing the
pseudo-effective reduction of $D$. Write $\Delta_0$ for the strict transform of $\Delta$ plus the reduced $\pi_0$-exceptional divisor and put $E_0=K_{W_0}+\Delta_0-\pi_0^*L$. There is a class $\gamma\in\Mov_1(W_0)$ with
\[
    f_*\gamma=0,\qquad \pi_{0*}\gamma=\alpha,\qquad E_0\cdot\gamma=0.
\]
\end{proposition}

\begin{proof}
Take the tower $\pi_m:W_m\to X$ of \cref{lem:approx} dominating $W_0$, with nef approximants $P_m$, and denote the maps to $W_0$ by $h_m$. Take simultaneous log resolutions so that the strict transform of $\Delta$ together with the reduced exceptional divisor has simple normal crossings on each $W_m$. Denote its boundary, with the coefficients of $\Delta$ retained on its strict transform, by $\Delta_m$. The numerical pullback $\alpha_m=\pi_m^!\alpha$ is movable and has zero intersection with every exceptional divisor. Indeed $\alpha$ is movable by extremality and $L\cdot\alpha<0$, and pushforward preserves pseudo-effective divisor classes. In particular
\[
  0\le S_m\cdot\alpha_m\le1 \qquad\text{by \cref{lem:approx}.}
\]

Fix $0<\epsilon<1$ sufficiently small and put
\[
L_\epsilon=K_X+(1-\epsilon)\Delta,\qquad T_m=K_{W_m}+(1-\epsilon)\Delta_m=\pi_m^*L_\epsilon+E_{m,\epsilon}.
\]
These pairs on $W_m$ are all $\epsilon$-lc, since they are log smooth and every boundary coefficient is at most $1-\epsilon$. Set $E_{0,\epsilon}=K_{W_0}+(1-\epsilon)\Delta_0-\pi_0^*L_\epsilon$. The discrepancy formula gives
\[
 E_{m,\epsilon}\ge h_m^*E_{0,\epsilon}\ge0.
\]
For the first inequality, the difference is exceptional over $W_0$ and its coefficients are the log discrepancies of the exceptional primes over $(W_0,(1-\epsilon)\Delta_0)$ minus $\epsilon$. The second inequality follows by the same calculation over $X$. In particular $T_m\cdot\alpha_m=L_\epsilon\cdot\alpha<0$.

Use the exact klt cone formula on $W_m$ to express $\alpha_m$ as the sum of a class in $\overline{\operatorname{NE}}_1(W_m)_{T_m\ge0}$ and finitely many terms $t[C]$ from Mori fibre rays. On each such ray choose Das's movable integral representative, so that
$0<-T_m\cdot C\le G(\dim X,\epsilon)$ \cite[Proposition~7.1]{Das20}. Retain the terms satisfying
\[
  -T_m\cdot C\ \ge\ \epsilon H\cdot \pi_{m*}C.
\]
Let $u_m$ be their sum and let $v_m=\alpha_m-u_m$. Each discarded term, and the original nonnegative term, has nonnegative degree against $T_m+\epsilon \pi_m^*H$. Consequently
\begin{equation}\label{eq:selection-budget}
(L_\epsilon+\epsilon H)\cdot(\alpha-\pi_{m*}u_m)\ge E_{m,\epsilon}\cdot u_m\ge E_{0,\epsilon}\cdot h_{m*}u_m.
\end{equation}
Here $E_{m,\epsilon}\cdot\alpha_m=0$ and $u_m$ is movable. All retained curves satisfy $1\le H\cdot \pi_{m*}C\le G(\dim X,\epsilon)/\epsilon$. They also satisfy $K_{W_m}\cdot C\le T_m\cdot C<0$ by movability.

Write $u_m=u_m^{\rm vert}+u_m^{\rm hor}$ by grouping the retained terms according to $f\circ h_m$. At this stage $\epsilon$ is fixed, so $M_\epsilon=G(\dim X,\epsilon)/\epsilon$ is a fixed bound. \cref{lem:vertical} gives $d_m(M_\epsilon)\to\infty$. Since $S_m$ is nef and $\alpha_m-u_m^{\rm hor}$ is in the closed effective cone, the horizontal coefficients satisfy, for all sufficiently large $m$,
\[
\sum_{C\text{ horizontal}} t_C\le \frac{S_m\cdot\alpha_m}{d_m(M_\epsilon)}\le \frac{1}{d_m(M_\epsilon)}.
\]
If the eligible set is empty the sum is zero. It follows that
\[
0\le \pi_0^*H\cdot (h_m)_*u_m^{\rm hor}\le M_\epsilon/d_m(M_\epsilon)\longrightarrow0.
\]
The big divisor $\pi_0^*H$ lies in the interior of the pseudo-effective divisor cone. Its unit slice on the dual cone $\Mov_1(W_0)$ is compact. Thus $(h_m)_*u_m^{\rm hor}\to0$ in $N_1(W_0)$. The full classes $(h_m)_*u_m$ are movable and have $\pi_0^*H$-degree at most one. Any cluster point $\gamma_\epsilon$ is consequently movable and vertical. Apply \eqref{eq:selection-budget} to the full $u_m$ along this subsequence. Since $\alpha-\pi_{m*}u_m=\pi_{m*}v_m$ is pseudo-effective, the limiting class $b_\epsilon=\alpha-\pi_{0*}\gamma_\epsilon$ satisfies
\[
b_\epsilon\in\overline{\operatorname{NE}}_1(X),\qquad(L_\epsilon+\epsilon H)\cdot b_\epsilon\ge E_{0,\epsilon}\cdot\gamma_\epsilon\ge0.
\]

Finally $\pi_0^*H\cdot\gamma_\epsilon\le1$. The same compactness on $W_0$ gives a subsequential limit $\gamma$ as $\epsilon\downarrow0$. Then $f_*\gamma=0$ and $b=\alpha-\pi_{0*}\gamma$ is pseudo-effective with $L\cdot b\ge E_0\cdot\gamma\ge0$. Since $\pi_{0*}\gamma$ is movable, $\alpha=b+\pi_{0*}\gamma$ is a decomposition inside
$\mathcal C_L(X)$. Extremality and $L\cdot\alpha<0$ give $b=0$. The same inequality now yields $E_0\cdot\gamma=0$. This last limit takes place entirely on $W_0$. No estimate on $d_m(M_\epsilon)$ uniform in $\epsilon$ is used.
\end{proof}

\subsection{Induction and descent}
\begin{proof}[Proof of \cref{thm:exact-dlt}]
First prove the assertion for log smooth pairs. In dimension one a negative adjoint degree on a smooth projective curve forces the curve to be rational. In positive dimension let $R$ be a negative extremal ray. When $\eta>0$, \cref{lem:div} gives a rational representative. When $\eta=0$, take the nonzero orthogonal movable divisor class $D$ of that lemma. On a smooth model realizing the pseudo-effective reduction of $D$, the positive part $P_\sigma(\pi_0^*D)$ is numerically trivial on a very general fibre \cite{LehmannReduction14}. Its pushforward to $X$ equals $D$ numerically: the negative part of $\pi_0^*D$ is exceptional, since $D$ is movable in codimension one. Consequently the reduction base has positive dimension. This remains true after simultaneous resolution of the boundary.

By \cref{prop:lift}, there is a nonzero movable class $\gamma$ on $W_0$ with $f_*\gamma=0$, $(\pi_0)_*\gamma=\alpha$, and $E_0\cdot\gamma=0$. Choose a very general point of $Z$ over which $f$ is smooth and every stratum of $\Delta_0$ meets the fibre transversely, and write $i:F\hookrightarrow W_0$ for this fibre. Peternell's fibre theorem in \cite[Theorem~4.5(iii)]{FulgerLehmann16} gives $u\in\Mov_1(F)$ with $i_*u=\gamma$. In particular $F$ has positive dimension, and it has dimension less than $X$. The boundary $\Delta_F=\Delta_0|_F$ is simple normal crossings by generic smoothness of all strata, and $T_0|_F=K_F+\Delta_F$ for $T_0=K_{W_0}+\Delta_0$.

The induction hypothesis, including its finite-sum formula, applies to $(F,\Delta_F)$. Since $u\in\Mov_1(F)\subset\NM(F)$, write
\[
u=b_F+\sum_i s_i[C_i],\qquad b_F\in\NE(F),\quad (K_F+\Delta_F)\cdot b_F\ge0,\quad s_i>0,
\]
where each $[C_i]$ spans a negative extremal ray of $\mathcal C_{K_F+\Delta_F}(F)$. For each such ray choose $C_i$ by the ray assertion in the induction hypothesis, so that its normalization $g_i:\PP^1\to F$ is free. The free maps deform in covering families on
$F$ and retain their numerical classes.

Along $g_i$, the tangent sequence of $f$ has trivial quotient:
\[
0\longrightarrow g_i^*T_F\longrightarrow g_i^*T_{W_0}\longrightarrow \cO_{\PP^1}^{\oplus\dim Z}\longrightarrow0.
\]
Since $g_i^*T_F$ is globally generated, $H^1(\PP^1,g_i^*T_F)=0$, so this sequence splits and $g_i$ is free on $W_0$. Its full deformation family dominates $W_0$ and preserves its numerical class. Thus each $C_i$ moves on $W_0$ with the same pushed-forward numerical
class. Take general members for all the numerical computations.

Let $b_0=i_*b_F$. Then $T_0\cdot b_0\ge0$. The effective divisor $E_0$ has nonnegative degree on each of the movable rational terms, so $E_0\cdot b_0\le E_0\cdot\gamma=0$. Therefore
\[
 L\cdot(\pi_0)_*b_0=(T_0-E_0)\cdot b_0\ge0.
\]
Pushing the decomposition of $\gamma$ to $X$ gives $\alpha$ as a finite sum of a class in $B_L(X)$ and movable rational classes. A general member of a dominating family meets the birational isomorphism locus, so its image is again a movable rational curve. The degree and extremality of $\alpha$ force at least one such curve to span $R$. A general normalization map in the resulting covering family on the smooth $X$ is free by the same evaluation argument. This completes the dimension induction for log smooth pairs.

Now suppose $(X,\Delta)$ is $\Q$-factorial dlt. Fix a negative extremal ray generated by $\alpha$. Choose a projective resolution $\pi:W\to X$ witnessing the dlt condition \cite[Definition~10.1]{Fujino11}. Its exceptional locus is a simple normal crossings divisor and, for $\Theta=\pi_*^{-1}\Delta+\sum_j E_j$,
\[
T=K_W+\Theta=\pi^*L+E,\qquad E=\sum_j A(E_j,X,\Delta)E_j,\qquad A(E_j,X,\Delta)>0.
\]
The numerical pullback $\pi^!\alpha$ is defined since $X$ is $\Q$-factorial. It is movable and annihilates $E$. By the log smooth case write $\pi^!\alpha=b+\sum_i t_i[C_i]$, where $b\in\NE(W)$, $T\cdot b\ge0$, $t_i>0$, and each $C_i$ is movable rational with $T\cdot C_i<0$. Each $E\cdot C_i\ge0$, so $E\cdot b\le0$ and $L\cdot\pi_*b\ge0$. The rational terms push to movable classes. Since $\alpha$ is negative extremal we get $\pi_*b=0$, whence
\[
0\le T\cdot b=E\cdot b=-\sum_i t_i E\cdot C_i\le0.
\]
Consequently $E\cdot C_i=0$ for every $i$. As all $A(E_j,X,\Delta)$ are positive, a general member of each covering family is disjoint from every $E_j$. It lies entirely in the isomorphism locus. Its general normalization is free on the smooth $W$, as in the fibre argument. Thus it gives a free rational curve in $X_{\rm sm}$; its nonzero class spans $R$ by the same extremality. For a free map the evaluation morphism from $\mathbb P^1$ times its deformation space to the smooth locus is smooth along the map. After shrinking the deformation space this holds at all points of the source curve. The inverse image of a fixed codimension-two subset has codimension at least two, so its projection misses a general parameter. This gives the stated avoidance property.

Finally take a small projective crepant $\Q$-factorialization $q:(\widehat X,\widehat\Delta)\to(X,\Delta)$ of the given dlt pair \cite[Theorem~10.4]{Fujino11}. A movable lift of $\alpha$ exists by \cite[Lemma~2.1]{Lehmann12}. Apply the preceding finite-sum formula to the lift. Crepancy sends its nonnegative term into $B_L(X)$. Each free rational term can be deformed away from the exceptional locus of $q$, which has codimension at least two, and away from the strict transform of any prescribed codimension-two closed subset. Its image is a free rational curve in $X_{\rm sm}$. Extremality and $L\cdot\alpha<0$ give a nonzero such term on $R$. The compact-slice argument in \cref{subsec:mmp-rays} now gives the asserted finite sums and countability.
\end{proof}

\end{document}